\documentclass[a4paper, 12pt]{article}
\usepackage{amsmath}
\usepackage{amsfonts}
\usepackage{amssymb}
\usepackage{amscd}
\usepackage{amsthm}

\newcommand{\G}{{\mathbb G}}

\newcommand{\A}{{\mathbb A}}
\renewcommand{\P}{{\mathbb P}}
\newcommand{\Z}{{\mathbb Z}}

\newcommand{\Q}{{\mathbb Q}}
\newcommand{\C}{{\mathbb C}}
\newcommand{\F}{{\mathcal F}}
\newcommand{\E}{{\mathcal E}}
\def\LL{\mathcal L}
\newcommand{\T}{{\mathcal T}}
\newcommand{\I}{{\mathcal I}}
\def\OO{\mathcal O}
\newcommand{\Tor}{{\mathcal Tor}}
\newcommand{\Ttor}{{\rm Tor}}
\newcommand{\Lef}{{\rm Lef}}

\newcommand{\m}{\mathfrak{m}}
\renewcommand{\a}{\mathfrak{a}}
\newcommand{\End}{{\rm End}}
\newcommand{\Id}{{\rm Id}}
\newcommand{\Tr}{{\rm Tr}}
\newcommand{\Hom}{{\rm Hom}}
\newcommand{\ch}{{\rm ch}}

\newcommand{\Sym}{{\rm Sym}}

\theoremstyle{theorem}
\newtheorem{theor}{Theorem}
\theoremstyle{theorem}
\newtheorem{prop}{Proposition}
\theoremstyle{theorem}
\newtheorem{corol}{Corollary}
\theoremstyle{theorem}
\newtheorem{lemma}{Lemma}

\theoremstyle{definition}
\newtheorem{defin}{Definition}

\theoremstyle{remark}
\newtheorem{rmk}{Remark}
\theoremstyle{remark}
\newtheorem{examp}{Example}

\title{Adelic Lefschetz formula for the action of a one-dimensional 
torus}
\author{S. O. Gorchinskiy, A. N. Parshin}
\date{}

\begin{document}
\maketitle

\section{Introduction}

There exists a well-known Lefschetz formula for the number of fixed 
points in algebraic topology: if $f:X\to X$ is an endomorphism of a 
compact oriented $n$-dimensional manifold $X$ with a finite number 
$N_f$ of fixed points and the determinant $\det(1-df)$ is positive 
at all the fixed points then
$$
N_f=\sum_{i=0}^n(-1)^i\Tr(f^*|_{H^i(X,\Q)}).
$$

In algebraic geometry, there exist cohomologies of {\it coherent 
sheaves} on algebraic varieties. So the following natural question 
arises: what is {\it the Lefschetz number} 
$\Lef(X,\F,f)=\sum_{i=0}^n(-1)^i\Tr(f^*|_{H^i(X,\F)})$ of the 
cohomologies of a coherent sheaf $\F$ on a projective variety $X$ 
equal to, if we are given the {\it lift} of the endomorphism $f$ up 
to the sheaf $\F$? The answer on this question is given by a 
coherent (or holomorphic) Lefschetz formula \cite{AB, Don}.

Suppose that $\F$ is locally free, i.e. corresponds to a vector 
bundle $F$, that the set $Z$ of fixed points is finite, and that 
the intersection of the graph of the endomorphism with the diagonal 
is transversal at all the fixed points. Then the coherent Lefschetz 
formula has the form
$$
\Lef(X,\F,f)=\sum_{x \in
Z}\frac{\Tr(f|_{F_{x}})}
{\det(1-df|_{T^*_{x}})},
$$
where $F_{x}$ is the fiber of the bundle $F$ over point $x$ and 
$T_{x}$ is the tangent space to $X$ at point $x$.

In the present paper, we consider the case of {\it the action of a 
one-dimensional torus} $\G_m$ on $X$ instead of {\it one 
endomorphism} of the variety $X$. In this case we can consider 
every summand in the Lefschetz formula, corresponding to a fixed 
point $x\in Z$, as the trace of the action of the element 
$f\in\G_m$ on the infinite-dimensional space 
$\hat{\OO}_{X,x}\otimes F_x$, where
$\hat{\OO}_{X,x}$ is the local completed ring of the point $x$. 
This remark was done first in \cite{PS}.

On the other hand, there exists a theory of adelic complexes 
$\A_X(\F)^{\bullet}$ \cite{Bel, Hub} for cohomologies of coherent 
sheaves $\F$. It allows to construct in a canonical way flasque 
resolutions of coherent sheaves and, hence, to compute their 
cohomologies. Because of its functoriality, the group $\G_m$ acts 
on the adelic complex. For example, if $X$ is a curve and the sheaf 
$\F$ corresponds to the vector bundle $F$ then the complex
$\A_X(\F)^{\bullet}$ is equal to
$$
\F_{\eta}\oplus \prod_x
\hat{\OO}_{X,x}\otimes F_x \to {\prod_x}^{\prime}
K_x \otimes F_x
$$
Its cohomologies coincide with
$H^{\bullet}(X, \F)$. Here $K_x$ denotes the local field of the 
point $x$.

Thus, for the computation of the Lefschetz number it should be 
natural to define the {\it traces} of the action of $\G_m$ on each 
component of the adelic complex. This problem is rather non-trivial 
because of the fact that the considered linear spaces are 
infinite-dimensional over the field $k$. If it were solved then 
since the adelic complex contains $\prod_{x\in Z}
\hat{\OO}_{X,x}\otimes F_x$ as a summand, the proof of the coherent 
Lefschetz formula would be immediately reduced to the explanation 
of why the trace on the rest of the adelic complex is equal to 
zero. There are two heuristic considerations explaining why this 
should be true:
\begin{itemize}
\item The trace of the action of the group on any field is equal to 
zero
\item The trace of the action of the group on the ''permutational'' 
representation is equal to zero
\end{itemize}

The second statement is well known for finite groups in 
the finite-dimensional case while the first one can be reduced to 
the second one by the normal base theorem from Galois theory.

If we suppose that we have defined the trace of the action of the 
group on the components of the adelic complex that verifies these 
properties then it will be easy to see that it should be equal to 
zero on all the components exept those that enter the right hand 
side part of the Lefschetz formula. The analogous reasoning should 
also take part in the case of an arbitrary dimension.

A direct realization of this plan seems to be not so simple. 
In this paper a somewhat roundabout way is proposed and is realized 
for the case of a locally free sheaf on a non-singular projective 
variety. Besides, the Lefschets formula acquires the form
$$
\Tr(\G_m,H^{\bullet}(X,
\F))=\Tr(\G_m,\A_X^{fix}(\F)),
$$
where $\A_X^{fix}(\F)$ is the ''fixed'' part of the adelic complex, 
related to the set of fixed points (see definition 1). It coincides 
with the product $\prod_{x\in Z}
\hat{\OO}_{X,x}\otimes F_x$ that was considered above in the case 
of a curve.

The definition of the complex $\A_X^{fix}(\F)$ is given for any 
scheme $X$ so it is possible to suppose that the adelic Lefschetz 
formula is valid for any proper scheme over the field $k$.

\section{Endomorphisms and actions of algebraic groups}

Let $X$ be a nonsingular projective algebraic variety $X$ defined 
over the field $k$. Let us consider an endomorphism $f:X\to X$ and 
a coherent sheaf $\mathcal F$ on $X$ with a chosen morphism (in the 
category of coherent sheaves) $\alpha:f^*\F\to\F$. We will call 
this morphism {\it a lift} of $f$ up to $\F$. We obtain an action 
of $f$ on the cohomologies of the sheaf $\F$ by taking the 
composition of the canonical map $H^i(X,\F)\to H^i(X,f^*\F)$ with 
the induced map 
$\alpha _*:H^i(X,f^*\F)\to H^i(X,\F)$. 

\begin{examp}\label{tangentaction} If $\F$ is the cotangent bundle 
$\T^*_X$ on $X$ then we will have the canonical lift 
$\alpha:f^*\T^*_X\to\T^*_X$ that is dual to the differential 
$df:\T_X\to f^*\T_X$.
\end{examp}

\begin{examp}\label{O(N+1)}
Let $f:\P^N\to\P^N$ be a projective automorphism. There is a 
canonical lift on the sheaf $\OO_{\P^N}(N+1)$ since 
$\OO_{\P^N}(N+1)\cong\wedge^N\T_X$. Thus, for all $n\in\Z$ there 
exist lifts $\alpha_n$ on $\OO_{\P^N}(n(N+1))$ such that 
$\alpha_n\otimes\alpha_m=\alpha_{n+m}$. In particular, 
$\alpha_1\otimes\alpha_{-1}=\Id$, where $\Id$ is the identity lift 
on the trivial bundle.
\end{examp}

For the situation described above we define the 
{\it Lefschetz number} by
$$
\Lef(X,\F,\alpha)=\sum_{i\ge 0}(-1)^i \Tr(f|_{H^i(X,\F)}).
$$

If $f$ acts trivially on the variety $X$ and on the sheaf $\F$ as 
well then $\Lef(X,\F,f)=\chi(X,\F)$, and the Lefschetz number (as 
the element of the ground field $k$) is given by the Riemann-Roch 
theorem \cite{BS},\cite{SGA6}.

\begin{rmk} If we summarize alternatively the Lefschetz numbers of 
the external powers of the cotangent bundle for a non-singular 
projective complex variety, then by Hodge decomposition we will get 
the ''classical'' Lefschetz number 
$\sum_{i\ge 0}(-1)^i \Tr(f|_{H^i(X,\C)})$.
\end{rmk}

In 1969 Donovan (see \cite{Don}) proposed to consider for a given 
variety $X$ and its endomorphism $f$ the category $Coh(X,f)$ whose 
objects are coherent sheaves $\F$ on $X$ together with lifts 
$\alpha$. The morphisms in this category are defined in a natural 
way:
$\Hom((\F_1,\alpha_1), (\F_2,\alpha_2))$ consists of such morphisms 
of coherent sheaves $\varphi:\F_1\to\F_2$ that the following 
diagram is commutative:
$$
\begin{array}{ccc}
f^*\F_1 & \stackrel{f^*\varphi}{\longrightarrow} & f^*\F_2 \\
\downarrow\lefteqn{\alpha_1}&&\downarrow\lefteqn{\alpha_2} \\
\F_1 & \stackrel{\varphi}{\longrightarrow} & \F_2 \\
\end{array}
$$
It is possible to consider an analogous category $LocFr(X,f)$ where 
coherent sheaves are replaced by locally free ones.

\begin{prop}\label{resol}
Let $f$ be an automorphism and suppose that there exists such a 
closed embedding of our variety into the projective space 
$X\hookrightarrow \P^N$ that $f$ can be extended to the projective 
automorphism of the whole $\P^N$. Then
$$
K_0(Coh(X,f))=K_0(LocFr(X,f)).
$$
\end{prop}
\begin{proof} The idea of the proof is absolutely the same as for 
the usual $K_0(X)$. It follows from example \ref{O(N+1)} that there 
exist canonical lifts up to the sheaves $\OO_X(n(N+1))$. It is 
enough to show that for any object $(\F,\alpha)$ of the category 
$Coh(X,f)$ one can construct its finite resolution using the 
objects of the category $LocFr(X,f)$. To do this we twist $\F$
by $\OO_X(n(N+1))$ (in the category $Coh(X,f)$) such that 
$\F(n(N+1))$ would be generated by its global sections. Then we get 
a covering $p:\oplus_i\OO_X\cdot e_i\to\F(n(N+1))$ where $e_i$ 
denote the elements of some base of $H^0(X,\F(n(N+1)))$. The action 
of $\alpha$ on $H^0(X,\F(n(N+1)))$ defines the lift of $f$ up to a 
trivial bundle
$\oplus_i\OO_X\cdot e_i$. To get the resolution of the initial 
sheaf $\F$ we just twist it back and then repeat the described 
procedure several times for the kernels that will appear. We will 
stop in a finite number of steps since $X$ is non-singular.
\end{proof}

\begin{rmk} The condition of proposition \ref{resol} is equivalent 
to the existence of an isomorphic lift up to an ample sheaf, i.e. 
to the existence of an embedding $X\hookrightarrow\P^N=\P(V)$ with 
the extension of the action of $\G_m$ on $X$ to a {\it linear} 
action of $\G_m$ on $V$ (it is the consequence of example 
\ref{O(N+1)}).
\end{rmk}

If $f$ acts on $X$ identically then the condition of proposition 
\ref{resol} is obviously satisfied. If the sheaf $\F$ is locally 
free then $\alpha$ defines linear operators on each fiber of the 
corresponding bundle $F\to X$. Since the coefficients of the 
characteristic polynomials of the action on fibers are regular 
functions on all of $X$, then they are constant and the 
characterictic polynomials are the same for all fibers. After we 
have passed to the algebraic closure of the field $k$ (we denote it 
by the same letter) we can consider the decomposition of a given 
linear space together with a given operator into Jordan blocks. The 
decomposition of fibers will also be ''constant over the base'', 
i.e. the bundle $F$ decomposes into the sum over eigenvalues:
$$
F=\oplus_{\lambda}F_{\lambda},
$$
where $\lambda\in k$ are the roots of the characteristic polynomial 
of an arbitrary fiber and $F_{\lambda}$ are maximal subbundles such 
that the corresponding operator $\lambda\cdot Id-\alpha$ is 
nilpotent on them.

After we have applied the functor $K_0$, the Jordan blocks 
decompose in their order into the sum of one-dimensional spaces and 
we get the isomorphism of rings:
$$
\begin{array}{rcl}
K_0(Coh(X,id_X))=K_0(LocFr(X,id_X))& \cong
& K_0(X)\otimes_{\Z} \Z[k] \\
(F, \alpha) &\mapsto & \sum_{\lambda}F_{\lambda}\otimes\lambda,\\
\end{array}
$$
where $\Z[k]$ denotes the group ring of the multiplicative 
semi-group of the field $k$. We need to consider it for two 
reasons. First, the direct sum of two isomorphic bundles with 
different lifts is not equal in the $K_0$-group to the bundle 
itself with the sum of lifts, and, secondly, the tensor product of 
bundles together with lifts results in pairwise multiplication of 
their eigenvalues.

There exists a natural homomorphism from the ring $\Z[k]$ to the 
field $k$ defined by the following formula:
$$
\begin{array}{rcl}
w: \Z[k]&{\to}&k \\
\sum_{\lambda} n_{\lambda}\cdot[\lambda]&
{\mapsto}& \sum_{\lambda} n_{\lambda}\lambda.\\
\end{array}
$$
Taking the composition of the constructed above isomorphism with 
the map $\chi\otimes w$, where $\chi:K_0(X)\to\Z$ is the Euler 
characteristic of a sheaf, we obtain the map $\chi_L$:
$$
\begin{array}{rcl}
\chi_L: K_0(Coh(X,id_X))&{\to}& \Z\otimes_{\Z} k=k \\
(F,\alpha)&{\mapsto}&
\sum_{\lambda}\chi(F_{\lambda})\cdot\lambda.\\
\end{array}
$$
Obviously,
$$
\Lef(X,\F,id_X)=\chi_L(\F).
$$
Thus, we see that for the computation of the Lefschetz number 
instead of knowing the element of the ring $K_0(X,id_X)$ itself, it 
is enough to know its image in $K_0(X)\otimes_{\Z}k$ since the 
Lefschetz number is well defined for the elements of the latter 
ring.

Now we will formulate a very useful criterion of the invertibility 
of elements of the ring $K_0(X)\otimes_{\Z}k$. Recall that there 
exists a canonical decomposition 
$K_0(X)\cong \Z\oplus\widetilde{K_0(X)}$, where the $\Z$-component 
corresponds to the virtual dimension of elements of $K_0(X)$. 
Moreover, the ideal $\widetilde{K_0(X)}$ is nilpotent (see 
\cite{SGA6}). It follows that the element of $K_0(X)$ is invertible 
if and only if its virtual dimension is an invertible integer, i.e. 
is equal to $\pm 1$. Then it is clear that there is an analogous 
criterion for the ring $K_0(X)\otimes_{\Z}k$ where the 
$k$-component will be equal to the trace of the action on an 
arbitrary fiber of a vector bundle $\F$ with a lift $\alpha$. Thus, 
the element $(F,\alpha)$ is invertible if and only if its trace on 
any fiber is not equal to zero.
 


Now let us return to the case where $f$ is not identity. If $f$ has 
a finite order then the set of fixed points $Z$ is the disjoint 
union of non-singular subvarieties $Z_{\alpha}$. Example 
\ref{tangentaction} shows that there is a canonical action of $f$ 
on the conormal bundle $N^*_{Z_{\alpha}/X}$. Moreover, $df$ does 
not act identically on any vector of any fiber
$N^*_{Z_{\alpha}/X}|_x$, $x\in Z_{\alpha}$.

An elementary fact from linear algebra says that for a vector space 
$V$ and an operator $A$ acting on $V$ the following indentity is 
true:
$$
\det (1_V-A)=\sum_{i=0}^{\dim V}(-1)^i \Tr(A|_{\wedge^i V}).
$$
Consequently, our criterion of invertibility implies that the 
element 
$\sum_{i\ge 0}(-1)^i{\wedge}{^i}N^*_{Z_{\alpha}/X}$ is invertible 
in the ring $K_0(Z_{\alpha})\otimes_{\Z}W(k)$. 

It is possible to show that the element 
$\sum_{i\ge0}(-1)^i{\wedge}{^i}N^*_{Z_{\alpha}/X}$
is invertible in the ring $K_0(Z_{\alpha})\otimes_{\Z}k$. Donovan 
proved the following 

\begin{theor} For an automorphism $f$ of a finite order and its 
lift $\alpha$ up to a coherent sheaf $\F$, the following identity 
holds true:
$$
\Lef(X,\F,\alpha)=\sum_{\alpha}\chi_L(\F|_{Z_{\alpha}}\otimes
(\sum_{i\ge0}(-1)^i {\wedge}{^i}N^*_{Z_{\alpha}/X})^{-1}),
$$
where by the restriction $\F|_{Z_{\alpha}}$ we mean the inverse 
image in $K_0$ groups with the induced action, in particular it 
coincides with the usual restriction $\F|_{Z_{\alpha}}$ if the 
sheaf $\F$ is locally free.
\end{theor}

\begin{rmk} This theorem may be reformulated as follows:
$$
\Lef(X,\F,\alpha)=\Lef(Z,\F|_Z\otimes
NL_{Z/X}, \alpha|_Z\otimes df_Z),
$$
where $NL_{Z/X}=(\sum_{i\ge 0}
(-1)^i{\wedge}{^i}N^*_{Z_{\alpha}/X})^{-1}$, and $df_Z$ denotes the 
canonical lift of $f$ up to this ''bundle''.

It is easy to see that if $\F$ is locally free, i.e. corresponds to 
a vector bundle $F$, and if $Z$ is a finite set then this formula 
will coincide with the Lefshetz formula that was mentioned in the 
introduction and was first proved by Atiah and Bott who used 
analytic methods \cite{AB}. The proof of Donovan is completely 
algebraic and does not depend on the properties of the field at 
all.
\end{rmk}

Now we consider an algebraic action of an algebraic group $G$ on a 
projective non-singular variety $X$ over the field $k$. Analogously 
to the case of the unique morphism, we consider a family of lifts 
corresponding to each element of the group $G$. The condition 
saying that it depends ''algebraically'' on the group may be stated 
as follows: there exists a morphism of sheaves $\alpha:m^*\F\to 
p^*\F$, where $m:X\times G\to X$ is the map that defines the action 
of the group, and $p:X\times G\to X$ is the projection on the 
second component. Thus, the case of the action of the group on the 
variety may be reduced, in some sense, to the case of the unique 
morphism. In particular, using it, we define the category 
$Coh(X,G)$ of the lifts of the action of the group up to $\F$. 
Morphisms in this category are such morphisms of coherent sheaves 
$\varphi:\F_1\to\F_2$ that the following diagram is commutative:
$$
\begin{array}{ccc}
m^*\F_1 & \stackrel{m^*\varphi}{\longrightarrow} & m^*\F_2 \\
\downarrow\lefteqn{\alpha_1}&&\downarrow\lefteqn{\alpha_2} \\
p^*\F_1 & \stackrel{p^*\varphi}{\longrightarrow} & p^*\F_2 \\
\end{array}
$$

We will also consider the category of lifts on only locally free 
sheaves $LocFr(X,G)$. If the variety $X$ is embedded in $\P^N$ and 
the action of $G$ on $X$ may be extended to the action on $\P^N$, 
then, like in proposition 2.1, we will get the isomorphisms of 
rings
$$
K_0(Coh(X,G))\cong K_0(LocFr(X,G)).
$$ 
Suppose that $G$ acts on the sheaf $\F$, i.e. that for any two 
elements $f,g\in G$ we have $\alpha_{fg}=\alpha_f\circ\alpha_g$. 
Such a sheaf is called {\it $G$-linearized}.
This may be reformulated in terms of the morphism 
$\alpha:m^*\F\to p^*\F$, giving a 1-cocycle condition on it (for 
details see \cite{Mum}). For a locally free sheaf this means in 
terms of the bundle that the group $G$ acts on its space $F$ 
(compatible with the action on the base $X$). For every $i$ we get 
a representation $\rho:G\to\End_k(H^i(X,\F))$.

Suppose that, like in the case of the unique morphism, the action 
of $G$ on $X$ is trivial and that the sheaf is locally free. Then 
we obtain a representation of the group in each fiber of the 
corresponding bundle. In fact, it is just a formalization of the 
notion of an "algebraic family of representations of the group 
$G$". It is possible to consider the decomposition of the vector 
bundle into the sum 
$\oplus_{\{\lambda(g)\}}F_{\lambda(g)}\otimes\lambda(g)\in 
K_0(X)\otimes_{\Z}k$ over eigenvalues as in the previous section 
for every element $g\in G$, by passing to the algebraic closure of 
the field $k$. It would be very nice if this decomposition depended 
''algebraically'' on $g$, i.e. if we could get some element of the 
ring $K_0(X)\otimes_{\Z}k[G]$, that would reduce to
$\oplus_{\{\lambda(g)\}}F_{\lambda(g)}\otimes\lambda(g)$ for every 
$g\in G$. Unfortunately, in the case of an arbitrary group G this 
decomposition doesn't exist. However, if $G$ is an algebraic torus 
$T$, then $F$ may be decomposed into the sum over eigenvalues for 
all elements of the torus at the same time. Namely,
$$
\oplus_{\chi}F_{\chi}\otimes\chi\in K_0(X)\otimes k[T],
$$
where $\chi$ runs over characters of the torus $T$ and $F_{\chi}$ 
is the biggest subbundle inside $F$ such that the action of any 
element $t\in T$ on it is the multiplication by $\chi(t)$. This 
follows from the fact that for any finite-dimensional 
representation of the torus there exists such an element $t_0\in T$ 
that all the characters of this representation are uniquely defined 
by their values in $t_0$ --- this is the algebraic analogue of the 
dense path on the real torus (cf. \cite{Car}). Hence, the 
decomposition of the bundle for the element $t_0$ will define the 
decomposition for the whole group.
For the bundle $F$ together with the action of the torus $T$ on it 
we will denote this decomposition by $\psi(F)$.

\begin{rmk}\label{invert}
It follows from the formulated criterion for the inversibility of 
elements of the ring $K_0(X)$ that the class $\psi(F)$ is 
inversible if and only if its trace on an arbitrary fiber is not 
identically zero as a function on the torus $T$. 
\end{rmk}

Now let us return to the case of a nontrivial action of the torus 
on the variety. As in the case of a finite group the following 
esult holds:
 
\begin{prop} The set of fixed points $Z=X^T$ is the dijoint union 
of non-singular subvarieties $Z_{\alpha}$. 
\end{prop}

The proof is given in \cite{Bia}. Informally speaking, it comes 
from the fact that the tangent space to a point of $X^T$ coincides 
with the fiber of the component of $\T_X|_{X^T}$, on which the 
torus acts trivially. Its dimension, of course, is constant on each 
connected component of $X^T$.

\begin{rmk}\label{invertable}
From what was said before and from remark \ref{invert} it follows 
that the element 
$\psi(\sum_{i\ge 0}(-1)^i{\wedge}{^i}N^*_{Z_{\alpha}/Z^-
_{\alpha}})$ is invertible in the ring $K_0(X^T)\otimes_{\Z}k(T)$.
\end{rmk}

The action of the torus also has another propertiy that is very 
useful for our purposes.

\begin{prop}[see \cite{Mum}]\label{linearization} If the torus $T$ 
acts on the projective variety $X$ then there exists an embedding 
$X\hookrightarrow\P^N$ and the extension of the action up to $\P^N$ 
that is compatible with the embedding.
\end{prop}

\begin{corol}\label{equiv}
There is an isomorphism of rings for the action of the torus $T$ on 
a non-singular projective variety $X$:
$$
K_0(Coh(X,T))\cong K_0(LocFr(X,T)).
$$
\end{corol}

From now on we suppose that the torus $T$ is a one-dimensional 
torus $\G_m$. Since every orbit is a regular map from 
$\G_m=\P^1\backslash\{0,\infty\}$ to the projective variety $X$, it 
is possible to define the ''limits'' $0(x),\infty(x)\in X$ for 
every point $x\in X$. This allows us to define for every connected 
component $Z_{\alpha}$ the set of ''incoming'' points 
$Z^+_{\alpha}$ consisting of points $x\in X$ such that $0(x)\in 
Z_{\alpha}$. In a similar way we define the set of ''outgoing'' 
points 
$Z^-_{\alpha}$ as the set of points $x\in X$ such that 
$\infty(x)\in Z_{\alpha}$.

\begin{examp} \label{simpleaction}
Let $X=\P^1$ and the element $\lambda\in\G_m$ act by the formula 
$\lambda(x_0:x_1)=(\lambda x_0:x_1)$. Then $Z=\{0\}\cup\{\infty\}$, 
$Z^+_{0}=\A^1=\{x_0\ne 0\}$, and $Z^+_{\infty}=\{\infty\}$.
\end{examp}

This example illustrates the general situation.

\begin{prop}[see \cite{Bia}] The subsets $Z^+_{\alpha}$ and 
$Z^-_{\alpha}$ are locally closed inside $X$ and have the structure 
of a non-singular algebraic variety. Besides, 
$X=\cup_{\alpha}Z^+_{\alpha}=\cup_{\alpha}Z^-_{\alpha}$ and 
$Z^+_{\alpha}\cap Z^-_{\alpha}=Z_{\alpha}$ for every $\alpha$.
\end{prop}

\begin{corol} For every $\alpha$ there exists a neighbourhood $U$, 
containing $Z_{\alpha}$, such that $Z^+_{\alpha}$ and $Z^-
_{\alpha}$ are closed inside $U$. For the ideal sheaves of 
corresponding subvarieties in $U$ the identity 
$\I_{Z_{\alpha}}=\I_{Z^+_{\alpha}}+\I_{Z^-_{\alpha}}$ holds true.
\end{corol}

Consider the action of $\G_m$ on the restriction to $Z_{\alpha}$ of 
the tangent bundle $T_X$ of $X$. The non-zero characters of the 
torus are divided into positive and negative ones, i.e. into the 
ones that may be continiously extended to $0$, and the ones that 
may be continiously extended to $\infty$, respectively. Then using 
a natural notation the decomposition of $T_X|_{Z_{\alpha}}$ over 
the characters may be written in the following way:
$$ 
T_X|_{Z_{\alpha}}=(T_X|_{Z_{\alpha}})_0\oplus(T_X|_{Z_{\alpha}})_+
\oplus(T_X|_{Z_{\alpha}})_-.
$$
The following result easily comes from our definitions.

\begin{prop} The following equality of bundles holds: 
$(T_X|_{Z_{\alpha}})_0=T_{Z_{\alpha}}$,
$(T_X|_{Z_{\alpha}})_+=N_{Z/Z^+_{\alpha}}$,
$(T_X|_{Z_{\alpha}})_-=N_{Z/Z^-_{\alpha}}$,
where $N_{Z/Z^{\pm}_{\alpha}}$ denote the normal bundles on $Z$ to 
$Z^{\pm}_{\alpha}$.
\end{prop} 

Let us also remark that $N_{Z/X}=N_{Z/Z^+_{\alpha}}\oplus
N_{Z/Z^-_{\alpha}}$.

\section{Adelic construction}

For an algebraic variety $X$ over an arbitrary field $k$ one can 
define {\it the adelic complex} $\A_X^{\bullet}$ 
\cite{Bel, Hub, FP}. We recall that $\A_X^{\bullet}$ is the complex 
of abelian groups and its $i$-th component is equal to
$$
\A_X^i={\prod_{\{\eta_0,\ldots,\eta_i\}}}{^{\prime}}
\hat{\OO}_{\eta_0,\ldots,\eta_i },
$$ 
where the product $\prod^{\prime}$ is taken over all flags of 
length $i$, i.e. over such sets of schematic points 
$\{\eta_0,\ldots,\eta_i\}$ that $\eta_{j+1}\in\overline{\eta_j}$ 
and $\eta_{j+1}\ne\eta_j$ for every index $j$. The ring 
$\hat{\OO}_{\eta_0,\ldots,\eta_i }$ is defined as follows. Take the 
completed local ring $\hat{\OO}_{\eta_i}$ and complete it on the 
ideal corresponding to $\eta_{i-1}$ (this ideal may not remain 
prime). Using this construction we obtain the ring 
$\hat{\OO}_{\eta_{i-1},\eta_i}$. Then we repeat this procedure $i-
1$ times for other schematic points in the flag and obtain the ring 
$\hat{\OO}_{\eta_0,\ldots,\eta_i}$. 

The sign $'$ means that instead of considering not the hole direct 
product over flags we consider a subset inside it satisfying some 
special {\it adelic condition} \cite{Bel, Hub, FP}.
The differential $d:\A_X^i\to\A_X^{i+1}$ is defined by the usual 
formula:
$$
(da)_{\eta_0,\ldots,\eta_i}=\sum_{j=0}^i 
(-1)^j a_{\eta_0,\ldots,\eta_{j-1},\eta_{j+1},\ldots,\eta_i}.
$$
Moreover, it is possible to define the adelic complex 
$\A_X(\F)^{\bullet}$ of a (quasi)coherent sheaf $\F$ on $X$ in a 
similar way. To do this we use modules 
$\hat{\F}_{\eta_0,\ldots,\eta_i}$ over rings 
$\hat{\OO}_{\eta_0,\ldots,\eta_i }$ that are sequential completions 
of the module $\F_{\eta_i}$. Many notions of algebraic geometry can 
be reformulated in terms of higher-dimensional adeles. For 
instance, Serre duality, Chern classes, Chow groups and 
intersection numbers \cite{FP}. In particular, the following 
general statement is true.

\begin{theor}\label{cohom} There exists a functorial isomorphism 
$H^*(X,\F)\cong H^*(\A_X(\F)^{\bullet}).$
\end{theor}

Consider the action of a one-dimensional torus $\G_m$ on a 
projective non-singular variety $X$ and a $\G_m$-linearized sheaf 
$\F$ on $X$.

\begin{defin}
We denote by $\A_X^{fix}(\F)^{\bullet}$ the quotient of the adelic 
complex whose $i$-th component is equal to
$$
\A_X^{fix}(\F)^{i}={\prod_{\{\eta_0,\ldots,\eta_i\}}}{^{\prime}}
\hat{\OO}_{\eta_0,\ldots,\eta_i },
$$
where the product is taken over all flags that lie inside the set 
of fixed points $X^{\G_m}$.
\end{defin}

\begin{examp} \label{simplefixed}
Suppose that we have finitely many fixed points. Then the complex 
$\A_X^{fix}(\F)^{\bullet}$ consists of only one component of degree 
zero. It will be equal to the sum of completed local rings of fixed 
points: 
$$
\A_X^{fix}(\F)^{\bullet}=\A_X^{fix}(\F)^0=
\oplus_{x\in Z}\hat{\OO}_x.
$$
\end{examp}

Since the ideal sheaves $\I_{Z^+_{\alpha}}$ and 
$\I_{Z^-_{\alpha}}$ are well-defined in some open subset $U$ 
containing $Z$, we get the bifiltration on the complex 
$\A_X^{fix}(\F)^{\bullet}$: 
$$
\A_{p,q}^{\bullet}=\A_X^{fix}(\F(\I_{Z^+_{\alpha}})^p
(\I_{Z^-_{\alpha}})^q)^{\bullet}.
$$
The quotients of this bifiltration may be found explicitly as shown 
below.

\begin{prop}\label{filtr} Suppose that $\F$ is locally free. Then
$$
\A_{p,q}^{\bullet}/(\A_{p+1,q}^{\bullet}+\A_{p,q+1}^{\bullet})\cong 
\A_Z((\F|_{Z}\otimes
\Sym^p(N^*_{Z/Z^+_{\alpha}})\otimes
\Sym^q(N^*_{Z/Z^-_{\alpha}}))^{\bullet}.
$$ 
\end{prop}
\begin{proof} In fact, the bifiltration is defined not only for 
adelic complex but also for the corresponding sheaf itself (on the 
open subset $U$). On the other hand the taking of the quotient of a 
finitely generated module $M$ over the ideal $I$ commutes with the 
completion:
$$
(\lim_{\longleftarrow}M/\m^iM)/(I\widehat{R_{\m}})=
\lim_{\longleftarrow}M/((\m^i+I)M).
$$
Consequently, it is enough to prove the required identity for the 
sheaf. 

Since $\F$ is locally free we have
$$
\F\I_{Z^+_{\alpha}}^p\I_{Z^-_{\alpha}}^q=
\F\otimes\I_{Z^+_{\alpha}}^p\I_{Z^-_{\alpha}}^q.
$$ 
Because $Z^+_{\alpha}$, $Z^-_{\alpha}$, and their ''infinitesimal 
neighbourhoods'' intersect transversely, 
we get $\Tor^{\OO_X}_i(\OO_X/\I^p_{Z^+_{\alpha}},\OO_X/\I^q_{Z^-
_{\alpha}})=0$ for $i>0$. Using the long exact sequence for 
$\Tor$-groups, it is easy to see that 
$\Tor^{\OO_X}_i(\I^p_{Z^+_{\alpha}},\I^q_{Z^-_{\alpha}})=0$ for 
$i>0$. From this it follows that
$$
\I^p_{Z^+_{\alpha}}\I^q_{Z^-_{\alpha}}\cong
\I^p_{Z^+_{\alpha}}\otimes\I^q_{Z^-_{\alpha}}.
$$
Finally, 
$$
\F\I_{Z^+_{\alpha}}^p\I_{Z^-_{\alpha}}^q/
(\F\I_{Z^+_{\alpha}}^{p+1}\I_{Z^-_{\alpha}}^q+
\F\I_{Z^+_{\alpha}}^p\I_{Z^-_{\alpha}}^{q+1})=
$$
$$
=\F\I_{Z^+_{\alpha}}^p\I_{Z^-_{\alpha}}^q\otimes
\OO_U/(\I_{Z^+_{\alpha}}+\I_{Z^-_{\alpha}})\cong
\F|_{Z_{\alpha}}\otimes\Sym^p(N^*_{Z/Z^+_{\alpha}}|_{Z_{\alpha}})
\otimes\Sym^q(N^*_{Z/Z^-_{\alpha}}|_{Z_{\alpha}}).
$$
Here we use the fact that if the quotient $R/I$ of the local 
regular ring $R$ remains regular, then 
$I^p/I^{p+1}\cong\Sym_{R/I}^p(I/I^2)$. The proof is finished.
\end{proof}

Let us remark that the torus $\G_m$ acts on all the complexes 
defined above: $\A_X(\F)^{\bullet}$, 
$\A_X^{fix}(\F)^{\bullet}$, $\A_{p,q}^{\bullet}$ and 
$\A_Z((\F|_{Z}\otimes
\Sym^p(N^*_{Z/Z^+_{\alpha}})\otimes
\Sym^q(N^*_{Z/Z^-_{\alpha}}))^{\bullet}$.

Let $V^{\bullet}$ be a bounded complex of vector spaces over $k$ 
with finite-dimensional cohomologies. Let us define the trace 
$\Tr(\G_m,V^{\bullet})$ of the action of $\G_m$ on $V^{\bullet}$ to 
be a regular function on $\G_m$ that is equal to the sum
$$
\Tr(\G_m,V^{\bullet})=
\sum_{i\in\Z}(-1)^i\Tr(\G_m|_{H^i(V^{\bullet})}).
$$
The problem is that this number is not well-defined by this formula 
if the cohomologies of the complex $V^{\bullet}$ are 
infinite-dimensional. To define the trace in this situation we 
embed the field of rational functions on $\G_m$ in the completed 
local field of some point lying on the compactification of the 
torus $\P^1=\G_m\cup\{0\}\cup\{\infty\}$. Now we can sum series of 
rational functions on the torus and define the trace as the sum of 
the traces on the quotients of some suitable filtration. In this 
manner one can extend the notion of the trace to the case of 
complexes with infinite-dimesional cohomologies. However, in this 
case the trace may become a rational function on the torus.

\begin{examp}\label{simpletrace} 
Let the complex consist of only one component --- a complete 
one-dimensional local ring $\OO$ with a local parameter $t$, and 
with the action of the torus $\lambda:t\mapsto\lambda t$ where 
$\lambda\in\G_m$. We have the filtration on $\OO$ by the degrees of 
the maximal ideal $\OO\supset t\OO\supset t^2\OO\supset\ldots$ 
while the torus acts on the $i$-th quotient of the filtration 
$t^i\OO/t^{i+1}\OO$ by multiplying it by $\lambda^i$. Therefore the 
trace of the action of $\G_m$ on $\OO$ should be equal to the 
series
$\sum_{i\ge 0}{\lambda}^i=\frac{1}{1-\lambda}$ that converges in 
the completed local ring $\hat{\OO}_{\P^1,0}\supset k[\G_m]$. 
\end{examp}

\begin{examp}
We join example \ref{simplefixed} with example \ref{simpletrace} 
and see that for the action from example \ref{simpleaction} the 
trace of the action of $\G_m$ on 
$\A_{\P^1}^{fix}(\OO_{\P^1})^{\bullet}$ is equal to 
$$
\frac{1}{1-\lambda}+\frac{1}{1-\lambda^{-1}}=1.
$$
On the other hand the Lefschetz number $\Lef(\P^1,\OO_X,\G_m)$ is 
also equal to 1: the sheaf $\OO_X$ has only one nontrivial 
cohomology $H^0$, which consists of constants on which the action 
of $\G_m$ is the identity.
\end{examp}

This example leads us to the following theorem.

\begin{theor}\label{biseries} Let the sheaf $\F$ be locally free. 
Then
\begin{enumerate}
\item For every natural $p$ the series 
$
\sum_{q\ge 0}\Tr(\G_m,
\A_{p,q}^{\bullet}/(\A_{p+1,q}^{\bullet}+\A_{p,q+1}^{\bullet})) 
$
converges in the complete ring $\hat{\OO}_{\P^1,\infty}$ to the 
rational function $\Tr_p\in k(\G_m)$; 

\item The series $\sum_{p\ge 0}\Tr_p$ converges in 
$\hat{\OO}_{\P^1,0}$ to the regular function $\Tr\in k[\G_m]$; 

\item The following equality is true:
$$
\Tr=\Lef(Z,\F|_Z\otimes NL_{Z/X},\G_m). 
$$
\end{enumerate}
\end{theor} 

\begin{proof} From the very beginning we pass to the algebraic 
closure of the field $k$. By proposition \ref{filtr} and theorem 
\ref{cohom} we have the identity
$$
\Tr_{p,q}=\Tr(\G_m,
\A_{p,q}^{\bullet}/(\A_{p+1,q}^{\bullet}+\A_{p,q+1}^{\bullet}))=
\Lef(Z, \F|_{Z}\otimes\Sym^p(N^*_{Z/Z^+_{\alpha}})\otimes
\Sym^q(N^*_{Z/Z^-_{\alpha}}),\G_m).
$$
Thus,
$$
\Tr_{p,q}=\sum_{\alpha}\chi_L(\F|_{Z_{\alpha}}\otimes
\Sym^p(N^*_{Z_{\alpha}/Z^+_{\alpha}})\otimes
\Sym^q(N^*_{Z_{\alpha}/Z^-_{\alpha}})).
$$

Consider the ring $K_0(X)\otimes\hat{\OO}_{\P^1,\infty}$. We define 
the dicrete valuation $\|\cdot\|$ on this ring by the formula 
$$
\|\sum_i F_i\otimes\lambda^{n_i}\|=\rho^{-\min_i\{n_i\}},
$$
where $\rho$ is an arbitrary real number bigger than 1. It is easy 
to check that $\|\cdot\|$ satisfies all the conditions of a 
discrete valuation and defines the structure of a metric space on 
the ring $K_0(X)\otimes\hat{\OO}_{\P^1,\infty}$. Let us denote its 
completion by $\hat{R}$. There exists a continious map 
$f:\hat{R}\to\hat{\OO}_{\P^1,\infty}$, that sends every element 
$\a\in K_0(X)\otimes\hat{\OO}_{\P^1,\infty}$ to a function in 
$\hat{\OO}_{\P^1,\infty}$ by the formula 
$$
f:\a\mapsto
\sum_{\alpha}\chi_L(\F|_{Z_{\alpha}}\otimes\Sym^p(N^*_{Z_{\alpha}/Z
^+_{\alpha}})
\otimes\a).
$$
Now let us prove the following statement.

\begin{lemma}\label{series}
The following identity is true in the topological ring $\hat{R}$:
$$
\sum_{q\ge 0}\Sym^q F\otimes t^q=(\sum_{i\ge 0}(-1)^i\wedge^i 
F\otimes t^i)^{-1}
$$
for any vector bundle $F$, where $t\in\hat{\OO}_{\P^1,\infty}$ is 
any element that has a positive valuation at $\infty$.
\end{lemma}
\begin{proof}
First, this series converges in $\hat{R}$ by the definition of a 
discrete valuation on $K_0(X)\otimes\hat{\OO}_{\P^1,\infty}$ and by 
the Cauchy criterion. 
Secondly, for any natural $k\ge 1$ there exists a well-known exact 
sequence 
$$
0\to\wedge^k F\to\wedge^{k-1}F\otimes F\to\wedge^{k-2}F\otimes 
\Sym^2F\to\ldots\to F\otimes\Sym^{k-1}F\to \Sym^k F\to 0
$$
This implies that
$$
(\sum_{q\ge 0}\Sym^q F\otimes t^q)\cdot(\sum_{i\ge 0}(-1)^i\wedge^i 
F\otimes t^i)=1.
$$
\end{proof}
Now let us recall that $\G_m$ acts on 
$N^*_{Z_{\alpha}/Z^-_{\alpha}}$ only with negative characters. 
Therefore we may decompose $N^*_{Z_{\alpha}/Z^-_{\alpha}}$ into a 
sum over characters 
$$
N^*_{Z_{\alpha}/N^-_{\alpha}}=
\oplus_{\chi_-}(N^*_{Z_{\alpha}/Z^-_{\alpha}})_{\chi_-}
$$
and apply lemma \ref{series} to every summand. 
Besides, since
$$
\sum_{q\ge 0}\psi(\Sym^q N^*_{Z_{\alpha}/Z^-_{\alpha}})=
\prod_{\chi_-}\sum_{q\ge 0}\Sym^q(N^*_{Z_{\alpha}/Z^-
_{\alpha}})_{\chi_-}\otimes\chi_-^q
$$
and
$$
\sum_{i\ge 0}(-1)^i\psi(\wedge^i N^*_{Z_{\alpha}/Z^-_{\alpha}})=
\prod_{\chi_-}\sum_{i\ge 0}(-1)^i\wedge^i (N^*_{Z_{\alpha}/Z^-
_{\alpha}})_{\chi_-}\otimes\chi_-^i,
$$
we obtain 
$$
\Tr_p=\sum_{q\ge 0}f(\psi(\Sym^q N^*_{Z_{\alpha}/Z^-_{\alpha}}))
=f(\sum_{q\ge 0}\psi(\Sym^q N^*_{Z_{\alpha}/Z^-_{\alpha}}))=
$$
$$
=f(\psi(\sum_{i\ge 0}(-1)^i
{\wedge}{^i}N^*_{Z_{\alpha}/Z^-_{\alpha}})^{-1}))=
\sum_{\alpha}\chi_L(\F|_{Z_{\alpha}}\otimes
\Sym^p(N^*_{Z/Z^+_{\alpha}})\otimes
(\sum_{i\ge 0}(-1)^i
{\wedge}{^i}N^*_{Z_{\alpha}/Z^-_{\alpha}})^{-1}),
$$
In particular, by remark \ref{invertable}, $\Tr_p$ belongs to the 
field $k(\G_m)$.

The same reasoning shows that
$$
\sum_{p\ge 0}\Tr_p=\sum_{\alpha}\chi_L(\F|_{Z_{\alpha}}\otimes
(\sum_{i\ge 0}(-1)^i{\wedge}{^i}N^*_{Z_{\alpha}/Z^+_{\alpha}})^{-
1}\otimes
(\sum_{i\ge 0}(-1)^i{\wedge}{^i}N^*_{Z_{\alpha}/Z^-_{\alpha}})^{-
1})=
$$
$$
=\sum_{\alpha}\chi_L(\F|_{Z_{\alpha}})\otimes
(\sum_{i\ge 0}(-1)^i{\wedge}{^i}N^*_{Z_{\alpha}/X})^{-1})=
\Lef(Z,\F|_Z\otimes NL_{Z/X},\G_m).
$$
Thus, the required indentity is proved.
\end{proof}

Now it is natural to give

\begin{defin}
$$
\Tr(\G_m,\A^{fix}_X(\F))=\Lef(Z,\F|_Z\otimes
NL_{Z/X},\G_m).
$$
\end{defin}

\begin{rmk}\label{reformul} The Lefschetz formula, written for the 
action of $\G_m$, has in this notation a very compact form:
$$
\Tr(\G_m,\A_X(\F))=\Tr(\G_m,\A_X^{fix}(\F)).
$$
In what follows we will prove this formula.
\end{rmk} 

\begin{rmk}\label{coherent}
We could consider an arbitrary coherent sheaf $\F$ on $X$ instead 
of a locally free one, construct the corresponding complex 
$\A_X^{fix}(\F)^{\bullet}$ and the same bifiltration on it 
$\A^{\bullet}_{p,q}=\A_X^{fix}(\F\I^p_{Z^+}\I^q_{Z^-})^{\bullet}$. 
It is easy to see that in this case as well the quotient of 
complexes 
$\A_{p,q}^{\bullet}/(\A_{p,q+1}^{\bullet}+\A_{p+1,q}^{\bullet})$ is 
equal to $\A_Z(\F\I^p_{Z^+}\I^q_{Z^-}/        
(\F\I^{p+1}_{Z^+}\I^q_{Z^-}+
\F\I^p_{Z^+}\I^{q+1}_{Z^-}))^{\bullet}$. So we could consider a 
biseries 
$$
\sum_p(\sum_q\Tr(\G_m,
\A_{p,q}^{\bullet}/(\A_{p,q+1}^{\bullet}+\A_{p+1,q}^{\bullet})))
$$
and suppose that it converges to the function
$\Lef(Z,\F|_Z\otimes NL_{Z/X},\G_m)$ where by the restriction 
$\F|_Z$ we mean the inverse image $K_0$ groups, i.e. this is the 
''usual'' restriction of some flat resolution of $\F$ on $X$. The 
convergence is supposed for one index after another in the rings 
$\OO_{\P^1,0}$ and $\OO_{\P^1,\infty}$ respectively (as in theorem 
\ref{biseries}). In this case the formula from remark 
\ref{reformul} would make sense for a coherent sheaf $\F$ as well.
\end{rmk}

\section{The exactness of the trace on bifiltrated complexes}

The main statement of this section is the following

\begin{theor}\label{exact}
Consider a $\G_m$-equivariant embedding $i:Y\hookrightarrow X$ of 
non-singular projective varieties and a locally free 
$\G_m$-equivariant sheaf $\F$ on $Y$. Let us construct a locally 
free resolution consisting of $\G_m$-equivariant sheaves to the 
coherent sheaf $i_*\F$ on $X$:
$$
0\to{\mathcal P}_n\to\ldots\to{\mathcal P}_0\to i_*\F\to 0.
$$
Then 
$$
\sum_{i=0}^n(-1)^i\Tr(\G_m,\A^{fix}_X({\mathcal P}_i)^{\bullet})=
\Tr(\G_m,\A_Y^{fix}(\F)^{\bullet}).
$$
\end{theor}

\begin{rmk}\label{intersection}
In the natural notation the indentities 
$Z_Y^{(\pm)}=Y\cap Z_X^{(\pm)}$ are true. These identities are also 
true for tangent spaces (this follows from the explicit description 
of tangent spaces to the varieties $Z^{(\pm)}$ in terms of the 
characters of the torus). Hence, the images of the ideal sheaves 
$\I_{Z^{(\pm)}_X}$ under the natural map $\OO_X\to \OO_Y$ are equal 
to $\I_{Z^{(\pm)}_Y}$.
\end{rmk}

\begin{proof}[Proof of theorem \ref{exact}]
Consider an exact sequence
$$
0\to\A_X^{fix}({\mathcal P}_n)^{\bullet}\to\ldots
\to\A_X^{fix}({\mathcal P}_0)^{\bullet}\to
\A_X^{fix}(\F)^{\bullet}\to 0.
$$
We have to prove the exactness of the trace, defined as the double 
limit over bifiltrations.

For this purpose we consider the quotient of complexes
$$
C_{\bullet}:0\to
\A_X^{fix}({\mathcal P}_n/\I_{Z^+}^p+\I_{Z^-}^q)^{\bullet}\to
\ldots\to
\A_X^{fix}({\mathcal P}_0/\I_{Z^+}^p+\I_{Z^-}^q)^{\bullet}\to
\A_X^{fix}(\F/\I_{Z^+}^p+\I_{Z^-}^q)^{\bullet}\to 0,
$$
where $Z=Z_X=X^{\G_m}$ and $\I_{Z^{\pm}}$ denote the ideal sheaves 
of incoming and outgoing components inside $X$ that are defined in 
some neighbourhood of $Z$. This sequence of complexes is not exact 
but has finite-dimensional cohomologies. Let $\Tr^{\prime}_{p,q}$ 
denote the alternated trace of the action of the torus $\G_m$ on 
this cohomologies. Remark \ref{intersection} implies that it is 
enough to prove that for every natural $p\ge 0$ there exists a 
limit $\Tr^{\prime}_p=\lim_{q\to\infty}{\Tr^{\prime}_{p,q}}$ in the 
local field $K_{\P^1,\infty}$ which takes a value in $k(\G_m)$, and 
that there exists a limit $\lim_{p\to\infty}{\Tr^{\prime}_p}$ in 
the local field $K_{\P^1,0}$ that is equal to zero.

Since the trace of the action of $\G_m$ on the bicomplex 
$C_{\bullet}$ is just the ''ordinary'' trace of the action on 
finite-dimensional cohomologies, we may consider the ''horizonal'' 
cohomologies $C_{\bullet}$ and compute the trace of the action of 
$\G_m$ on its cohomologies. Explicitly, ''horizontal'' cohomologies 
are equal to
$$
\A_X^{fix}(\Tor^{\OO_X}_i(i_*\F,\OO_X/(\I_{Z^+}^p+
\I_{Z^-}^q))^{\bullet},
$$
where $i\ge 1$. From the properties of the functor $\Tor$ it 
immediately follows that 
$\Tor^{\OO_X}_i(i_*\F,\OO_X/(\I_{Z^+}^p+\I_{Z^-}^q)$ is supported 
inside $Z_Y=Y^{\G_m}$. Thus, we see that
$$
\Tr^{\prime}_{p,q}=\sum_{i=1}^n(-1)^{i+1}\chi_L(\LL^{p,q}_i),
$$
where $\LL^{p,q}_i$ denote the sheaves 
$\Tor^{\OO_X}_i(i_*\F,\OO_X/(\I_{Z^+}^p+\I_{Z^-}^q))$.
Therefore, analogously to the proof of theorem \ref{biseries}, we 
need to prove that for every $i\ge 1$ there exists a limit 
$L^p_i=\lim_{q\to\infty}\ch_L(\LL^{p,q}_i)$ in the topological ring 
$K_0(X)\otimes\hat{\OO}_{\P^1,\infty}$ which takes a value in 
$K_0(X)\otimes k(\G_m)$ and also that there exists a limit
$\lim_{q\to\infty}L^p_i$ in the topological ring 
$K_0(X)\otimes\hat{\OO}_{\P^1,0}$, which is equal to 0. 

To begin with, we consider the fibers of the sheaves $\LL^{p,q}_i$ 
over an arbitrary point $y\in Z_Y$. Since $\F$ is locally free, we 
have an equality of $\G_m$-modules
$$
(\LL^{p,q}_i)_y=\Ttor^{\OO_{X,y}}_i(\F_y,\OO_{X,y}/(I^p_{Z^+} 
+I^q_{Z^-}))=
\F|_y\otimes 
\Ttor^{\OO_{X,y}}_i(\OO_{Y,y},\OO_{X,y}/(I^p_{Z^+}+I^q_{Z^-})),
$$ 
where $\F|_y$ denotes the fiber over $y$ of the vector bundle on 
$Y$ that corresponds to the locally free sheaf $\F$.

Let us make the explicit calculations of the second factor. 
Consider the following local rings: $R_+=\hat{\OO}_{Z^+,y}$, $R_-
=\hat{\OO}_{Z^-,y}$, $R_0=\hat{\OO}_{Z,y}$, $R=R_+\otimes_k R_-
\otimes_k R_0$, the maximal ideals $I_+\subset R_+$, 
$I_-\subset R_-$, and also the ideals of the subvariety $Y$ in the 
corresponding local rings
$J_+\subset R_+$, $J_-\subset R_-$ and $J_0\subset R_0$. It is 
clear that $\hat{R}=\hat{\OO}_{X,y}$, and that the completion of 
the ring $R_+/J_+\otimes_k R_-/J_-\otimes_k R_0/J_0$ coincides with 
$\hat{\OO}_{Y,y}$. Since the completion is an exact functor it 
commutes with the application of $\Ttor$. Let us also remark that 
the completion of a finite-dimensional vector space is the space 
itself. These two facts imply that
$$
\Ttor^{\OO_{X,y}}_i(\OO_{Y,y},\OO_{X,y}/(I^p_{Z^+}+I^q_{Z^-}))=
\Ttor^{\hat{\OO}_{X,y}}_i(\hat{\OO}_{Y,y},\hat{\OO}_{X,y}/(I^p_{Z^+
}+
I^q_{Z^-}))=
$$
$$
=\Ttor^R_i
(R_+/J_+\otimes_k R_-/J_-\otimes_k R_0/J_0,R_+/I_+^p\otimes_k R_-
/I_-^q \otimes_k R_0)
$$
All the algebras over the field are flat over it 
(even infinite-dimensional). So we have the identity
$$
\Ttor^R_{\bullet}
(R_+/J_+\otimes_k R_-/J_-\otimes_k R_0/J_0,R_+/I_+^p\otimes_k R_-
/I_-^q\otimes_k R_0)=
$$
$$
=\Ttor^{R_+}_{\bullet}
(R_+/J_+,R_+/I_+^p)\otimes_k
\Ttor^{R_-}_{\bullet}
(R_-/J_-,R_-/I_-^q)\otimes_k R_0/J_0.
$$
In particular, we obtain that $\LL^{p,q}_i$ are locally free 
sheaves on $Z_Y$ (a module over a local ring is free if and only if 
its completion is free because the completion is an exact functor).

Having considered minimal resolutions of modules 
$R_{\pm}/I_{\pm}^{p(q)}$ over local rings $R_{\pm}$, we see that 
the valuations at $0$ and $\infty$ of the characters of the 
represantation of $\G_m$ in 
$\Ttor^{R_{\pm}}_i(R_{\pm}/J_{\pm},R_{\pm}/I_{\pm}^p))$ are not 
less than $p$ and $q$ respectively. Thus, for every fixed $p$ and 
for sufficiently large $q$ the locally free sheaves $\LL^{p,q}_i$ 
can be decomposed into a sum over the characters of a ''rather 
large'' order at $\infty$ and a sum over all the remaining 
characters:
$$
\LL^{p,q}_i=(\LL^{p,q}_i)_1\oplus(\LL^{p,q}_i)_2.
$$
At the level of completed fibers it corresponds to the 
decomposition
$$
(\hat{\LL}^{p,q}_i)_y=
\oplus_{k>0} (\F|_y \otimes_k
\Ttor^{R_+}_{i-k}(R_+/J_+,R_+/I_+^p)\otimes_k 
\Ttor^{R_-}_k(R_-/J_-,R_-/I_-^q)\otimes_k R_0/J_0)\oplus
$$
$$
\oplus \F|_y \otimes_k
\Ttor^{R_+}_i(R_+/J_+,R_+/I_+^p)\otimes_k 
R_-/J_-\otimes_k R_-/I_-^q\otimes_k R_0/J_0.
$$
Thus, similarly to the proof of theorem \ref{biseries}, it follows 
from lemma \ref{series} that the limit 
$\lim_{q\to\infty}\ch_L((\LL^{p,q}_i)_1)=0$ is equal to zero in the 
ring $K_0(X)\otimes\hat{\OO}_{\P^1,\infty}$, and that there exists 
a limit $\lim_{q\to\infty}\ch_L((\LL^{p,q}_i)_2)$ that belongs to 
$K_0(X)\otimes k(\G_m)$. Moreover, it tends to zero when 
$p\to\infty$ in the ring $K_0(X)\otimes\hat{\OO}_{\P^1,0}$.
This finishes the proof of theorem \ref{exact}.
\end{proof}

\section{Main theorem and its proof}

Now we have developed our machinery sufficiently to prove the main 
statement of this paper.

\begin{theor}\label{main}
For every locally free $\G_m$-linearized sheaf $\F$ on $X$ the 
following equality is true:
$$
\Tr(\G_m,\A_X(\F)^{\bullet})=\Tr(\G_m,\A_X^{fix}(\F)^{\bullet}).
$$
\end{theor}

\begin{proof}
By proposition \ref{linearization} we may assume that we have a 
$\G_m$-equivariant embedding $i:X\hookrightarrow\P^N$. Besides, we 
may assume that the action of $\G_m$ on $\P^N$ is linear, i.e. that 
$\G_m$ acts on the vector space $V$ so that $\P^N=\P(V)$, or, 
equivalently, that the sheaf $\OO_{\P^N}(1)$ is 
$\G_m$-linearized. 

Consider the sheaf $i_*\F$ on $\P^N$, together with the lift of the 
action of $\G_m$ from $X$. The complexes $\A_X(\F)$ and 
$\A^{fix}_X(\F)$ do coincide with the complexes $\A_{\P^N}(i_*\F)$ 
and $\A^{fix}_{\P^N}(i_*\F)$ respectively, because, by definition, 
the fibers of the sheaf $i_*\F$ outside $X$ are equal to zero. 
Also, the bifiltration on the complex $\A^{fix}_X(\F)$ induced from 
$\P^N$ coincides with the bifiltration, induced from $X$: this 
comes from remark \ref{intersection}. 

We costruct a locally free resolution of the sheaf $i_*\F$ on 
$\P^N$:
$$
0\to{\mathcal P_n}\to\ldots\to{\mathcal P_0}\to i_*\F\to 0
$$
As it was shown in theorem \ref{exact}, the following equality 
holds true:
$$
\sum_{i=0}^n(-1)^i\Tr(\G_m,\A^{fix}_X({\mathcal P}_i))=
\Tr(\G_m,\A^{fix}_X(\F)).
$$

Thus, it is enough to prove the statement of the theorem only for 
locally free sheaves on $\P^N$. Moreover, the construction of the 
resolution in the proof of proposition \ref{resol} implies that it 
is enough to consider locally free sheaves of the form 
$\E=\oplus_{i=1}^n\OO_{\P^N}(l)$. Here the action of the torus is 
just a tensor product of some $n$-dimensional represantation $\rho$ 
of the torus $\G_m$ with the canonical action of $\G_m$ on 
$\OO_{\P^N}(l)$. It is clear that 
$$
\Tr(\G_m,\A_{\P^N}(\E))=\Lef(\P^N,\E,\G_m)=
\Lef(\P^N,\OO_{\P^n}(l),\G_m)\cdot\Tr(\rho),
$$
and
$$
\Tr(\G_m,\A^{fix}_{\P^N}(\E))=\Lef((Z,\E|_Z\otimes 
NL_{Z/\P^N},\G_m)=\Lef((Z,\OO_Z(l)\otimes NL_{Z/\P^N},\G_m) 
\cdot\Tr(\rho),
$$
where $Z=(\P^N)^{\G_m}$. Hence we have reduced the general case to 
the case of the sheaf $\OO_{\P^N}(l)$. Since all the linear 
representations of the torus are diagonalizable, the variety 
$Z_{\P^N}$ is just a disjoint union of projective subspaces 
$\coprod_{\alpha}\P^{N_{\alpha}}$ corresponding to the characters 
of the represantation of $\G_m$ in $V$. After all these remarks the 
case of the sheaf $\OO_{\P^N}(l)$ is reduced to a straightforward 
but rather dull calculations that are similar to those in \cite{BS} 
on page 597. We omit further details.
\end{proof}

\begin{corol}
For the action of $\G_m$ on the smooth projective variety $X$ and 
an arbitrary $\G_m$-linearized coherent sheaf $\F$ the following 
form of the Lefschetz trace formula in the above notation is true:
$$
\Lef(X,\F,\G_m)=\Lef(Z,\F|_Z\otimes NL_{Z/X},\G_m).
$$
\end{corol}
\begin{proof}
It is just a reformulation of theorem \ref{main} for locally free 
sheaves. To prove it for an arbitrary coherent sheaf we remark that 
both sides of the needed equality are exact on sheaves, and use a 
locally free resolution of a coherent sheaf.
\end{proof}

\begin{rmk}
Despite the fact that the Lefschetz formula is true for coherent 
sheaves it remains impossible to prove that coherent sheaves 
satisfy the condition of remark \ref{coherent}.
\end{rmk}

\section{Appendix: purely non-hyperbolic case}

In one special case it is possible to prove that coherent sheaves 
satisfy the condition of remark \ref{coherent}. In this case the 
bifiltration on the quotients of adelic complexes turns out to be 
just an ordinary linear filtration, and it becomes possible to 
prove the exactness of the trace using a more general fact about 
filtrated complexes.
 
\begin{defin}
We say that the action of a one-dimensional torus is {\it purely 
non-hyperbolic} if for any connected component $Z_{\alpha}$ of the 
set of fixed points $Z$ either $Z_{\alpha}^{+}$ or $Z_{\alpha}^{-}$ 
is empty.
\end{defin}

In this case the bifiltration on the complex 
$\A^{fix}_{\P^N}(\F)^{\bullet}$ is just a filtration by powers of 
either the ideal $\I_{Z^{+}_{\alpha}}$ or the ideal 
$\I_{Z^{-}_{\alpha}}$. 

\begin{defin}
Consider the representation $\rho$ of the torus $\G_m$ in an 
infinite-dimensional linear space $V$ over the ground field $k$. We 
say that {\it the order of $\rho$ is at least $k$ at point} 
$x\in\P^1=\overline{\G}_m$, if for any $\G_m$-invariant subspace 
$U\subset V$ and for any $\G_m$-equivariant finite-dimensional 
quotient $U\to W$ the oreders of all its characters as functions on 
$\G_m$ are at least $k$ at $x$.
\end{defin}

\begin{rmk}\label{defcohom}
Suppose that we have the action of the torus $\G_m$ on the complex 
$V^{\bullet}$ of linear spaces over $k$ and a $\G_m$-equivariant 
filtration $F^p V^{\bullet}$ on $V^{\bullet}$ so that the orders of 
the representations $F^p V^{\bullet}$ are at least $p$ at the point 
$x\in\P^1$, and such that the cohomologies of the quotients $F^p 
V^{\bullet}/F^{p+1} V^{\bullet}$ are finite-dimensional. It comes 
from our definitions that the orders of the traces of $\G_m$ on the 
cohomologies of $F^p V^{\bullet}/F^{p+1} V^{\bullet}$ 
are at least $p$. Then the series $\sum_p\Tr(\G_m,F^p 
V^{\bullet}/F^{p+1} 
V^{\bullet})$ converges in the completed local ring $\OO_{\P^1,x}$. 
The sum of this series will be denoted by 
$\Tr(\G_m,F^{\bullet}V^{\bullet})$
\end{rmk} 

\begin{examp}
Let $X$ be a non-singular projective variety with the action of 
$\G_m$ and $\F$ be a locally free $\G_m$-linearized sheaf on it. 
Then, up to an additive constant, the order of the action of $\G_m$ 
on each component of $\A_X^{fix}(\F\I^p_{Z^+}\I^q_{Z^-})^j$ of the 
complex 
$\A_X^{fix}(\F\I^p_{Z^+}\I^q_{Z^-})^{\bullet}$ is at least $p$ at 
$0$ and at least $q$ at 
$\infty$. This comes from the fact that on each component the 
bifiltration by the degrees of ideals $\I_{Z^{(\pm)}}$ induces the 
bifiltration on any finite-dimensional quotient of any subspace 
inside this component. Proposition \ref{filtr} implies that the 
representation of the torus on the quotients of the initial 
bifiltration decomposes into the finite direct sum of ''scalar'' 
representations, that is representations such that the torus acts 
on them by one character only. Besides, the orders of these 
characters are not less than $p$ at $0$ and not less than $q$ at 
$\infty$ up to an additive constant (this constant arises because 
of the characters of the action of $\G_m$ on $\F|_Z$). 

Since every coherent sheaf on $X$ may be $\G_m$-equivariantly 
covered by a locally free sheaf, the present statement remains 
valid for arbitrary coherent sheaves.
\end{examp}

\begin{corol}

Suppose that the action of the torus on $X$ is purely 
non-hyperbolic. Then for any coherent sheaf $\F$ on $X$ the series 
$$
\sum_p\Tr(\G_m,
\A_p^{\bullet}/\A_{p+1}^{\bullet}),
$$
corresponding to the filtration by the powers $\I_{Z^{\pm}}$,
converges in the local ring $\OO_{\P^1,0}$ or $\OO_{\P^1,\infty}$.
\end{corol}

Now we prove the following rather abstract statement. 

\begin{prop}\label{1dimexact}
Let
$$
0\to V^{\bullet}_0\to V^{\bullet}_1\to\ldots\to V^{\bullet}_n\to 0
$$
be an exact sequence of $\G_m$-complexes. Suppose we have a 
$\G_m$-equivariant filtration $F^p V^{\bullet}_i$ such that all 
together $F^p V^{\bullet}_i$ form a subcomplex in the initial 
sequence and satisfy the condition from remark \ref{defcohom} for 
some point $x\in\P^1$. Then 
$$
\sum_{i=0}^n(-1)^i\Tr(\G_m,F^{\bullet} V^{\bullet}_i)=0.
$$
\end{prop}
\begin{proof}
Consider the quotients of the complexes 
$V_i^{\bullet}/F^p V_i^{\bullet}$. Together they form a bicomplex. 
The spectral sequence, whose first term consists of 
finite-dimensional cohomologies of these quotients 
$V_i^{\bullet}/F^p V_i^{\bullet}$ with respect to a ''vertical'' 
differential, i.e. the differential from the adelic complex, tends 
to the cohomologies of the whole bicomplex. Let us remark that it 
does not matter on what finite-dimensional level of the spectral 
sequence we compute the trace of $\G_m$. Hence, the trace of the 
torus on the cohomologies of the bicomplex is equal to 
$$
\sum_{i=0}^n(-1)^i\sum_{j=0}^{p-1}
\Tr(\G_m,F^j V_i^{\bullet}/F^{j+1}V_i^{\bullet}).
$$
On the other hand we could compute the cohomologies of the 
bicomplex using the spectral sequence associated to the 
''horizontal'' differential. The ''horizontal'' cohomologies of 
$V_{\bullet}^j/F^p V_{\bullet}^j$ turn out to be isomorphic to 
''horizontal'' cohomologies of the subcomplex $F^p V_{\bullet}^j$ 
because the initial sequence of complexes was exact.  
Consequently, the orders of the action of the torus on them, i.e. 
the orders at $x$ of the action of the torus on the cohomologies of 
the bicomplex, are at least $p$ up to an additive constant. So, the 
order at $x$ of the (alternated) sum 
$$
\sum_{i=0}^n(-1)^i\sum_{j=0}^{p-1}
\Tr(\G_m,F^j V_i^{\bullet}/F^{j+1}V_i^{\bullet})
$$
tends to infinity when $p\to\infty$, and the required statement is 
proved.
\end{proof}

\begin{corol}\label{exactcorol}
Let the action of the torus on $X$ be purely non-hyperbolic and let 
$$
0\to\F_1\to\ldots\to\F_n\to 0
$$
be an exact $\G_m$-equivariant sequence of coherent sheaves on $X$. 
Then
$$
\sum_{i=0}^n(-1)^i\Tr(\G_m,\A^{fix}_X(\F_i))=0.
$$
\end{corol}

\end{document}